\documentclass[11pt, reqno]{amsart}
\usepackage{amsmath}
\usepackage{amssymb}

\font\tengothic=eufm10 \font\sevengothic=eufm7
\newfam\gothicfam
      \textfont\gothicfam=\tengothic
      \scriptfont\gothicfam=\sevengothic
\def\goth#1{{\fam\gothicfam #1}}

\numberwithin{equation}{section}
\newenvironment{proofc}{{\bf Proof of Theorem \ref{reg}.}}{\hspace*{\fill} $\Box$
\par\vspace{1ex}}

\begin{document}
\setlength{\baselineskip}{1.7em}
\newtheorem{thm}{\bf Theorem}[section]
\newtheorem{pro}[thm]{\bf Proposition}
\newtheorem{claim}[thm]{\bf Claim}
\newtheorem{lemma}[thm]{\bf Lemma}
\newtheorem{cor}[thm]{\bf Corollary}
\newtheorem{Remark}[thm]{Remark}
\newtheorem{Example}[thm]{Example}

\newcommand{\G}{{\mathcal G}}
\newcommand{\depth}{\mbox{\rm depth } }
\newcommand{\reg}{\mbox{reg } }
\newcommand{\ga}{\mathfrak{a}}
\newcommand{\gr}{\goth r}
\newcommand{\pp}{{\mathbb P}}
\newcommand{\zz}{{\mathbb Z}}
\newcommand{\nn}{{\mathbb N}}
\newcommand{\R}{{\mathcal R}}
\newcommand{\A}{{\mathcal A}}
\newcommand{\x}{{\mathbb X}}
\newcommand{\y}{{\mathbb Y}}
\newcommand{\ix}{I_{\mathbb X}}
\newcommand{\mi}{{\mathcal I}}
\newcommand{\B}{{\bf B}}
\newcommand{\bi}{{\bf I}}
\newcommand{\V}{{\bf V}}
\newcommand{\cv}{{\mathcal V}}
\newcommand{\bL}{{\bf L}}
\newcommand{\calo}{{\mathcal O}}
\newcommand{\cl}{{\mathcal L}}
\newcommand{\cm}{{\mathcal M}}
\newcommand{\cn}{{\mathcal N}}
\newcommand{\m}{-\!\!--\!\!\rightarrow}
\newcommand{\smap}{\rightarrow\!\!\!\!\!\rightarrow}
\newcommand{\sfrac}[2]{\frac{\displaystyle #1}{\displaystyle #2}}
\newcommand{\img}{\Lambda_{d+1}}
\newcommand{\under}[1]{\underline{#1}}
\newcommand{\ov}[1]{\overline{#1}}
\newcommand{\sheaf}[1]{{\mathcal #1}}
\newcommand{\lex}{{\le}_{\mbox{\scriptsize lex}}}

\title[Depth of the associated graded ring]{The depth of the associated graded ring of
ideals with any reduction number} \author{Ian Aberbach, Laura Ghezzi and Huy T\`ai H\`a}
\address{Department of Mathematics, University of Missouri,
Columbia MO 65211.}
\email{aberbach@math.missouri.edu \\ ghezzi@math.missouri.edu \\
tai@math.missouri.edu}
\begin{abstract}
Let $R$ be a local Cohen-Macaulay ring, let $I$ be an $R$-ideal,
and let $\G$ be the associated graded ring of $I$. We give an
estimate for the depth of $\G$ when $\G$ is not necessarily
Cohen-Macaulay. We assume that $I$ is either equimultiple, or has
analytic deviation one, but we do not have any restriction on the
reduction number. We also give a general estimate for the depth of
$\G$ involving the first $\gr+\ell$ powers of $I$, where $\gr$
denotes the Castelnuovo regularity of $\G$ and $\ell$ denotes the
analytic spread of $I$.

\vspace*{3ex}
\end{abstract}
\maketitle

\setcounter{section}{-1}
\section{Introduction}

Let $R$ be a Noetherian local ring with infinite residue field
$k$, and let $I$ be an $R$-ideal. The {\it Rees algebra}
${\mathcal R}=R[It] \cong \oplus_{i\geq 0} I^i$ and the {\it
associated graded ring} ${\mathcal G}=gr_I(R)=\mathcal R \otimes_R
R/I \cong \oplus_{i\geq 0} I^i/I^{i+1}$ are two graded algebras
that reflect various algebraic and geometric properties of the
ideal $I$. For example, Proj$(\mathcal R)$ is the blow-up of
Spec$(R)$ along $V(I)$ and Proj$({\mathcal G})$ corresponds to the
exceptional fiber of the blow-up. Many authors have extensively
studied the Cohen-Macaulay property of $\mathcal R$ and ${\mathcal
G}$. The most general results have been obtained by Johnson and
Ulrich \cite[3.1]{AN2} and by Goto,
 Nakamura and Nishida \cite[1.1]{Goto}.
 The goal of this paper is to estimate the depth of $\mathcal G$
and $\mathcal R$
 when these rings are not necessarily Cohen-Macaulay.
 We can focus on the study of $\depth
{\mathcal G}$, since if ${\mathcal G}$ is not Cohen-Macaulay, we
have that $\depth {\mathcal R}$ = $\depth {\mathcal G}+1$
\cite[3.10]{HM}. In order to state and motivate our results, we
first need to recall some definitions and background.
\newline A very useful tool in the study of blow-up rings is the
notion of reduction of an ideal, with the reduction number
measuring how closely the two ideals are related. This approach is
due to Northcott and Rees \cite{NR}. An ideal $J\subseteq I$ is
called a {\it reduction} of $I$ if the morphism
$R[Jt]\hookrightarrow R[It]$ is finite, or equivalently if
$I^{r+1}=JI^r$ for some $r\geq 0$. The least such $r$ is denoted
by $r_J(I)$. A reduction is {\it minimal} if it is minimal with
respect to inclusion, and the {\it reduction number} $r(I)$ is
defined as $\min\{r_J(I)\mid J\ {\rm a\ minimal\ reduction\ of\ }
I\}$. One of the big advantages of reductions is that they contain
a lot of information about the ideal $I$, but often require fewer
generators. More precisely, every minimal reduction of $I$ is
generated by $\ell$ elements, where $\ell=\ell(I)$ is the analytic
spread of $I$; i.e., the Krull dimension of the ring ${\mathcal
R}\otimes_R k \cong {\mathcal G} \otimes_R k$. The analytic spread
is at least the height $g$ of $I$, and at most the dimension of
$R$. The difference $\ell-g$ is the {\it analytic deviation} of
$I$. Ideals for which the analytic deviation is zero are said to
be {\it equimultiple}. For further details see
\cite{Wolmer}.\newline Cortadellas and Zarzuela came up with
formulas for $\depth \G$ in \cite{cz}, in the special cases of
ideals with analytic deviation at most one and reduction number at
most two. Ghezzi in \cite{laura} found a general estimate of
$\depth \G$ involving the depth of the powers of the ideal $I$ up
to the reduction number (see \cite[2.1]{laura} for the precise
statement). This theorem recovers the formulas of \cite{cz} and
generalizes the results of \cite{AN2} and \cite{Goto}. However, in
the set-up of \cite{laura} (as well as in \cite{AN2}, \cite{Goto},
and \cite{cz}), the reduction number is at most the ``expected''
one. Namely, the assumptions of \cite[2.1]{laura} imply that
$r(I)\le \ell-g+1$. The main goal of this paper is to find an
estimate of $\depth \G$ without any restriction on $r(I)$. In
Section 1 we treat the cases in which the ideal is either
equimultiple, or has analytic deviation one. We make an assumption
on $\depth_{\G_+} \G$, where $\G_+$ denotes the ideal of $\G$
generated by homogeneous elements of positive degree. We are now
ready to state our main results.

\noindent{\bf Theorem~\ref{equi}.} \ {\it Let $R$ be a local
Cohen-Macaulay ring with infinite residue field, and let $I$ be an
equimultiple ideal with height $g$ and reduction number $r$. Let
$t=\min \{\depth R/I^j - r +j | 1 \le j \le r \}$.
\newline
{\rm (1)} If $\depth_{\G_+} \G = g$, then $\depth \G \ge g +
\max\{ 0,t\}$.\newline {\rm (2)} If $\depth_{\G_+} \G = g-1$, then
$ \depth \G \ge g +\max\{ -1,t\}$.}

 In particular, if the reduction
number is small, we have a formula for $\depth \G$.

\noindent{\bf Corollary~\ref{r2}.} {\it Let $R$ be a local
Cohen-Macaulay ring with infinite residue field, and let $I$ be an
equimultiple ideal with height $g$ and reduction number two.
Assume that $\depth R/I^2< \depth R/I$.
 If either $\depth_{\G_+} \G = g$ or $\depth_{\G_+} \G = g-1$, then $ \depth \G = g
+\depth R/I^2$.}

\noindent{\bf Theorem~\ref{dev1}.} {\it Let $R$ be a local
Cohen-Macaulay ring with infinite residue field, and let $I$ be an
analytic deviation one ideal with height $g$ and reduction number
$r$. Assume that $I$ is generically a complete intersection, and
that $\depth_{\G_+} \G = g$. Let $t=\min \{\depth R/I^j - r +j | 1
\le j \le r \}$. Then, $ \depth \G \ge g + 1 + \max\{ -1,t\}.$}

The key fact in the proofs of Theorem \ref{equi} and of Theorem
\ref{dev1} is that we can reduce to the case where $\depth_{\G_+}
\G = 0$ and so the reduction of the ideal is principal. The case
of a reduction generated by two elements is more complicated (see
Proposition~\ref{equi2} for a special case).\newline In Section 2
we give a lower bound for $\depth \G$ in terms of the depth of the
first $\gr+\ell$ powers of the ideal $I$. Here $\gr$ denotes the
 Castelnuovo-Mumford regularity of the associated
graded ring of $I$. In general it is known that $\gr \ge r(I)$,
but our results of Section 2 are valid for ideals with any
reduction number (not just ideals with the expected reduction
number). We first recall the definition and some notation that we
will use throughout Section 2.\newline
 Let $S = \oplus_{n \ge 0} S_n$ be a finitely
generated standard graded ring over a Noetherian ring $S_0$. For
any graded $S$-module $M = \oplus_{n \ge 0} M_n$, we define
\[ a(M) := \left\{\begin{array}{ll} \max \{ n | M_n \not= 0 \} & \mbox{ if } M \not= 0, \\

-\infty & \mbox{ if } M = 0. \end{array} \right. \] Let
$S_+=\oplus_{n > 0} S_n$ be the ideal generated by the homogeneous
elements of positive degree of $S$. For $i \ge 0$, set
\[ a_i(S_+,S) := a(H^i_{S_+}(S)), \]
where $H^i_{S_+}(.)$ denotes the $i$th local cohomology functor
with respect to the ideal $S_+$. The {\it Castelnuovo-Mumford
regularity} of $S$ is defined as the number
\[ \reg S := \max \{ a_i(S_+,S) + i | i \ge 0 \}. \]
This is an important invariant of the graded ring $S$ (see for
instance \cite{trung} and the literature cited there).
\newline
The main result of Section 2 can be stated as follows.

\noindent{\bf Theorem~\ref{reg}.} {\it Let $R$ be a local
Cohen-Macaulay ring with infinite residue field, and let $I$ be an
$R$-ideal with analytic spread $\ell$. Let $\G $ be the associated
graded ring of $I$, and $\gr = {\rm reg}\  \G$. Then, $ \depth \G
\ge \min (\{ \depth R/I^j | 1 \le j \le \gr+1 \}\cup \{ \depth
R/I^j+j-\gr | 2+\gr \le j \le \ell+\gr \} ). $ }

The proof of Theorem~\ref{reg} uses the techniques of
\cite{laura}. The result is inspired by work of Trung
\cite{trung}, that shows that we can find a minimal reduction of
$I$ with ``good intersection properties" (see
Lemma~\ref{filter-regular}).

\section{Main Results}

The following theorem gives a lower bound of $\depth \G$ for
equimultiple ideals with any reduction number.

\begin{thm} \label{equi}
Let $R$ be a local Cohen-Macaulay ring with infinite residue
field, and let $I$ be an equimultiple ideal with height $g$ and
reduction number $r$. Let $t=\min \{\depth R/I^j - r +j | 1 \le j
\le r \}$.
\begin{enumerate}
\item
 If $\depth_{\G_+} \G = g$, then $\depth \G \ge g + \max\{
0,t\}$.
\item
 If $\depth_{\G_+} \G = g-1$, then $ \depth \G \ge g
+\max\{ -1,t\}$.
\end{enumerate}
\end{thm}
\begin{proof} We prove the results by induction on $\depth_{\G_+}
\G$.\newline
 {\bf (1)}  Suppose that $\depth_{\G_+} \G = 0$. Then
$I^{r+1} = 0$. Hence, $ \G = R/I \oplus I/I^2 \oplus \ldots \oplus
I^{r-1}/I^r \oplus I^r$, and $\depth \G = \min \{ \depth R/I,
\depth I/I^2, \ldots, \depth I^{r-1}/I^r, \depth I^r \}$.\newline
Since $\depth I^i/I^{i+1} \ge \min \{\depth R/I^{i+1}, \depth
R/I^i+1\}$ for each $1\le i\leq r-1$, and
$\depth I^r \ge \depth R/I^r$, we have that
$\depth \G \ge \min \{ \depth R/I^j | 1 \le j \le r \}\ge t$.
\newline Now assume that $\depth_{\G_+} \G > 0$. Let $x \in I$ be an element such that
$\ov{x} \in I/I^2$ is regular on $\G$. By \cite[2.7]{vv} $x$ is
regular on $R$ and $I^j \cap (x) = xI^{j-1}$ for every $j \ge 1$.
Let $\ov{R} = R/(x),\ \ov{I} = I/(x)$ and $\ov{\G} = \G/(\ov{x}) =
\G_{\ov{R}}(\ov{I})$. By the induction hypothesis, we have that $
\depth \ov{\G} \ge g-1 + \min\{ \depth \ov{R}/\ov{I}^j - r +j | 1
\le j \le r \}$, and so $ \depth \G \ge g + \min\{ \depth
\ov{R}/\ov{I}^j - r +j | 1 \le j \le r \}. $ For $2 \le j \le r$,
consider the exact sequence $$ 0 \rightarrow R/xI^{j-1}
\rightarrow R/I^j \oplus R/(x) \rightarrow \ov{R}/\ov{I}^j
\rightarrow 0.$$ It follows that $ \depth \ov{R}/\ov{I}^j  \ge
\min \{ \depth R/I^j, \depth R/I^{j-1}- 1\}$ for  $1 \le j \le r$.
Hence we have that $\depth \G \ge g+ t$.

{\bf (2)} Suppose that $\depth_{\G_+} \G = 0$. Let $J =(a)$ be a
minimal reduction of $I$ with $r_J(I)=r$. For every $j \ge r+1$ we
have an exact sequence
\[ 0 \rightarrow R/I^{j-1} \rightarrow R/I^j \rightarrow R/(a) \rightarrow 0. \]
Using induction on $j$ we see that $ \depth R/I^j \ge \depth
R/I^r$ for every $j \ge r$. Hence, $ \depth \G \ge \inf \{ \depth
R/I^j | j \ge 1 \}= \min \{ \depth R/I^j | 1 \le j \le r \}  \ge
t.$ We may assume that $t\ge 0$. If $t>0$, let $x_1,\dots, x_t\in
R $ be a regular sequence on $R$ and on $R/I^j$ for all
$j=1,\dots, r$. Write $\ov{R}=R/(x_1,\dots, x_t)$,
$\ov{I}=I\ov{R}$. Since $\depth \ov{R}/\ov{I}^j =\depth R/I^j-t$
for all $j=1,\dots, r$, we have that $\min \{ \depth
\ov{R}/\ov{I}^j-r+j| 1 \le j \le r \}=0$. Hence we can reduce the
problem to the case where $t=0$. Now, choose $x \in R$ such that
$x$ is regular on $R$ and on $R/I^j$ for all $j = 1,\ldots, r-1$.
Let $x^*$ and $a^*$ be the initial forms of $x$ and $a$ in $\G$
($x^*$ has degree 0 and $a^*$ has degree 1). We claim that
$x^*+a^*$ is regular on $\G$, which proves the assertion. If not,
there exists $v^* = v_0^* + \ldots+v_n^* \in \G$, with $v_i^* \in
I^i/I^{i+1}$ and at least one $v_i^* \not= 0$, such that
$(x^*+a^*)v^* = 0$ in $\G$. Suppose that $n \ge r-1$. Then,
$a^*v_n^* = 0$ implies that $av_n \in I^{n+2} = aI^{n+1}$, and so
$v_n \in I^{n+1}$ since $a$ is regular on $R$. Hence $v_n^* = 0$.
Let $v_k^*\neq 0$ be the lowest degree term of $v^*$. Then
$x^*v_k^* = 0$ implies that $xv_k \in I^{k+1}$. Since $k+1 \le
r-1$, $x$ is regular on $R/I^{k+1}$, and so $v_k \in I^{k+1}$,
i.e., $v_k^* = 0$, a contradiction. This finishes the proof of the
case $\depth_{\G_+} \G = 0$.\newline If $\depth_{\G_+} \G > 0$, we
can follow the same induction step as that of part (1) to prove
the theorem.
\end{proof}

The following remark gives an upper bound for $\depth \G$ in a
general context.

\begin{Remark}\cite[2.11]{laura}\label{upperbound} Let $R$ be a Noetherian
local ring, let $I$ be an $R$-ideal with analytic spread $\ell$.
Then $\depth \G \leq \inf \{\depth R/I^j\mid j\geq 1\}+\ell$.
\end{Remark}

The next corollary is a special case of Theorem \ref{equi}, for
reduction number two. Combining Theorem \ref{equi} with Remark
\ref{upperbound}, we have a formula for $\depth \G$.

\begin{cor} \label{r2}
Let $R$ be a local Cohen-Macaulay ring with infinite residue
field, and let $I$ be an equimultiple ideal with height $g$ and
reduction number two. Assume that $\depth R/I^2< \depth R/I$.
 If either $\depth_{\G_+} \G = g$ or $\depth_{\G_+} \G = g-1$, then $ \depth \G = g
+\depth R/I^2$.
\end{cor}

In the next example we compute the depth of the associated graded
ring.

\begin{Example} {\rm Let $R = k[[X, Y, T_1, \ldots, T_n]]/(X^3Y) = k[[x, y, t_1, \ldots,
t_n]]$, where $k$ is a field and $n \ge 2$. $R$ is Cohen-Macaulay
and $\dim R = n+1$. Let $I = (xy, t_1, \ldots, t_{n-1})$, and let
$J = (t_1, \ldots, t_{n-1})$. $I$ is equimultiple with $\mbox{ht }
I = n-1$ and reduction number 2. $J$ is a minimal reduction of
$I$. We have that $\depth_{\G_+} \G = n-1$ by \cite[2.7]{vv},
since $I^2 \cap J = JI$. Furthermore, $\depth R/I = 2$ and $\depth
R/I^2 = 1$. Corollary \ref{r2} implies that $\depth \G = n$.}
\end{Example}

In the next theorem we treat the case of analytic deviation one
ideals with any reduction number. We obtain a lower bound for
$\depth \G$ similar to that of Theorem \ref{equi}.

\begin{thm} \label{dev1}
Let $R$ be a local Cohen-Macaulay ring with infinite residue
field, and let $I$ be an analytic deviation one ideal with height
$g$ and reduction number $r$.
Assume that $r(I_{\wp}) < r$ for every prime $\wp$ containing $I$ with
$\mbox{\rm ht } \wp = g$,
and that $\depth_{\G_+} \G = g$. Let $t=\min \{\depth R/I^j - r + j | 1 \le j \le r \}$.
Then,
\[ \depth \G \ge g + 1 + \max\{ -1,t\}.\]
\end{thm}

\begin{proof} We prove the result by induction on $\depth_{\G_+} \G$.
Suppose that $\depth_{\G_+} \G = 0$, and let $J = (a)$ be a
minimal reduction of $I$ with $r_J(I)=r$. Since $(0:a) \cap I^r =
0$, by \cite[3.4]{cz} we have that $\depth R/I^{r+1} = \depth
R/I^r$, if $R/I^r$ is not Cohen-Macaulay, and that $\depth
R/I^{r+1} = \depth R/I^r - 1$, if $R/I^r$ is Cohen-Macaulay. In
particular $R/I^{r+1}$ is not Cohen-Macaulay, and so
\cite[3.4]{cz} implies that $\depth R/I^j = \depth R/I^{r+1}$ for
every $j \ge r+1$. Hence, if $R/I^r$ is not Cohen-Macaulay, then
$\depth \G \ge  \min  \{ \depth R/I^j | 1 \le j \le r \} \ge t$.
If $R/I^r$ is Cohen-Macaulay, then $\depth \G \ge  \min  (\{
\depth R/I^j | 1 \le j \le r-1 \} \cup \{\dim R-1\}) \ge t$, if
$r\ge 2$. If $r \le 1$, by \cite[3.1]{z} we still have
 that
$ \depth \G \ge t$. Now we proceed as in the proof of part (2) of Theorem
\ref{equi} to obtain that $\depth \G \ge 1 + t$. This finishes the
proof of the case $\depth_{\G_+} \G = 0$. \newline Suppose now that
$\depth_{\G_+} \G > 0$. We follow again the same induction step of
Theorem \ref{equi} to get the assertion.
\end{proof}

The following example is an application of Theorem \ref{dev1}.

\begin{Example}{\rm  Let $R = k[[X, Y, Z, W, T_1, \ldots, T_n]]/(X^4Y,ZW) =
 k[[x, y, z, w, t_1, \ldots,
t_n]]$, where $k$ is a field and $n \ge 2$. $R$ is Cohen-Macaulay
and $\dim R = n+2$. The ideal $I= (xy, z, t_1, \ldots, t_{n-1})$
has height $n-1$, analytic deviation 1, reduction number 3, and it
is generically a complete intersection. The ideal $J= (z, t_1,
\ldots, t_{n-1})$ is a minimal reduction of $I$. We have that
$\depth_{\G_+} \G = n-1$ by \cite[2.7]{vv}, since
$I^n\cap(t_1,\dots,t_{n-1})=(t_1,\dots,t_{n-1})I^{n-1}$ for every
$n\ge 1$ . Since $\depth R/I = 3$, $\depth R/I^2 = 2$, and $\depth
R/I^3 = 1$, Theorem \ref{dev1} and Remark \ref{upperbound} imply
that $\depth \G = n+1$.}
\end{Example}

We remark again that the key fact in the proofs of Theorem
\ref{equi} and of Theorem \ref{dev1} is that we can reduce to the
case where the reduction of the ideal is principal. Even when the
reduction is generated by a regular sequence of two elements the
situation is much more complicated. The next proposition treats a
special case.

\begin{pro} \label{equi2}
Let $R$ be a Noetherian local ring with infinite residue field,
and let $I$ be an equimultiple ideal of height two and reduction
number $r$. Let $J = (a_1, a_2)$ be a minimal reduction of $I$
such that $I^r : a_1 = I^r : a_2$. Then, $ \depth R/I^j \ge \min
\{ \depth R/I^r -1, \depth R/I^{r+1} \}$ for every $j \ge r+1$,
and
\[ \depth \G \ge \min ( \{ \depth R/I^j | 1 \le j \le r-1 \}\cup \{\depth R/I^r - 1,
\depth R/I^{r+1}\} ) \]
\end{pro}
\begin{proof} For $j \ge r+1$ we have that $I^j =J^{j-r}I^r$, and so we have the following
exact sequences
\begin{eqnarray}
0 \rightarrow \sfrac{R}{a_1J^{j-r-1}I^r \cap a_2^{j-r}I^r}
&\rightarrow& \sfrac{R}{a_1I^{j-1}} \oplus \sfrac{R}{a_2^{j-r}I^r}
\rightarrow R/I^j \rightarrow 0. \label{eq4}
\end{eqnarray}
Since $I^r : a_1 = I^r : a_2$, we have that $a_1J^{j-r-1}I^r \cap
a_2^{j-r}I^r = a_1a_2^{j-r}(I^r : J)$. The sequence (\ref{eq4})
for $j=r+1$ implies that
\[ \depth \sfrac{R}{I^r : J} \ge \min \{ \depth R/I^r, \depth R/I^{r+1} +1 \}. \]
Now we prove by induction on $j$ that $ \depth R/I^j \ge \min \{ \depth R/I^r-1, \depth
R/I^{r+1} \}$ for all $j \ge r+1$. The claim is clear for $j =
r+1$. Suppose that $j \ge r+2$. From the sequence (\ref{eq4}), we
have that
$\depth R/I^j \ge \min \{ \depth R/I^r-1, \depth R/I^{r+1}, \depth
R/I^{j-1} \} \ge \min \{ \depth R/I^r - 1, \depth R/I^{r+1} \}.$
The first assertion is proved. The second assertion follows from
the first one.
\end{proof}

\begin{Remark} {\rm  Let $R$ be a Noetherian local ring with infinite residue field, and
let $I$ be an analytic deviation one ideal of height one and reduction number
$r$. Let $J = (a_1, a_2)$ be a minimal reduction of $I$ such that
$(a_1: a_2^n) \cap I^r \subseteq (a_1)$ for every $n \ge 1$, and
$I^r : a_1 = I^r:a_2$. Then, by the proof of Proposition \ref{equi2},
we have that $\depth R/I^j \ge \min \{ \depth R/I^r -1, \depth R/I^{r+1} \}$
for every $j \ge r+1$, and $\depth \G \ge \min ( \{ \depth R/I^j | 1 \le j \le r-1 \} \cup
\{ \depth R/I^r - 1, \depth R/I^{r+1} \} )$.}
\end{Remark}

\section{Depth of $\G$ and its Castelnuovo-Mumford regularity}

Let $R$ be a local Cohen-Macaulay ring, let $I$ be an
$R$-ideal with analytic spread $\ell$, let $\G$ be the associated
graded ring of $I$, and let $\goth{r} = \reg\G$. The purpose of this
section is to give a lower bound for $\depth \G$ involving the
depth of the first $\gr+\ell$ powers of $I$. We first prove some
technical results (Lemmas~\ref{filter-regular}, \ref{l1}, and
\ref{l2}), that will play a crucial role for the proof of
 Theorem \ref{reg}.

\begin{lemma} \label{filter-regular}
Let $R$ be a local ring with infinite residue field, let $I$ be an
ideal of $R$, $\G=gr_I(R)$ and $\gr={\rm reg}\
\G$. For $a\in I$
let $a^*$ denote the image of $a$ in $[\G]_1$. Let $J$ be a
reduction of $I$. Then there exists a minimal basis $a_1, \ldots,
a_s$ of $J$ satisfying the following conditions:
\begin{eqnarray}
[(a_1, \ldots, a_i) : a_{i+1}] \cap I^j &=& (a_1,\ldots,
a_i)I^{j-1}, \ \forall \ 0 \le i \le s-1,\  j \ge \gr + 1.
\label{eqn2}
\end{eqnarray}
\begin{eqnarray}[(a_1^*, \ldots, a_i^*) : a_{i+1}^*]_j &=&
(a_1^*, \ldots, a_i^*)_j, \ \forall\   0 \le i \le s-1, \ j \ge
\gr + 1.\label{eqn21}\end{eqnarray}

\end{lemma}

\begin{proof}
By \cite[1.1]{trung} there exists a minimal basis $a_1, \ldots,
a_s$ of $J$ such that $[(a_1, \ldots, a_i): a_{i+1}] \cap
I^{\gr+1} = (a_1, \ldots, a_i)I^{\gr}$, whenever $0 \le i \le
s-1$. Thus, \cite[4.7]{trung} implies (\ref{eqn2}). By
\cite[4.8]{trung}, $a_1^*, \ldots, a_s^*$ is a filter-regular
sequence of $\G$. We have that $[(a_1^*, \ldots, a_i^*) :
a_{i+1}^*]_j = (a_1^*, \ldots, a_i^*)_j$ whenever $0 \le i \le
s-1$ and $j \ge a(a_1^*, \ldots, a_s^*) + 1$, where $a(a_1^*,
\ldots, a_s^*)$ is defined to be $\max\{ a[(a_1^*, \ldots, a_i^*)
: a_{i+1}^* / (a_1^*, \ldots, a_i^*)]\mid i=0,\dots, s-1\}.$ By
\cite[2.4]{trung} we have that $\gr \ge a(a_1^*, \ldots, a_s^*)$.
This implies (\ref{eqn21}).
\end{proof}

In particular, (\ref{eqn2}) implies that $[0:a_1] \cap I^j = 0 \
\forall \ j \ge \gr + 1$. We remark that since $k$ is infinite, we
can choose the basis $ a_1, \ldots, a_s$ such that each $a_i$
satisfies $ [0:a_i] \cap I^j = 0 \ \forall \ j \ge \gr+1. $

\begin{lemma} \label{l1}{\rm (see \cite[2.3]{laura})}
Let $R$ be a local Cohen-Macaulay ring of dimension $d$ with
infinite residue field, and let $I$ be an $R$-ideal. Let $J$ be a
reduction of $I$ with
 basis $a_1,\dots,a_s$ satisfying {\rm
(\ref{eqn2})}. Write $\ga_i = (a_1, \ldots, a_i)$ for all $i = 0,
\ldots, s$. Then, $ \depth R/\ga_i I^j \ge \min (\{ d-i\}\cup \{
\depth R/I^{j-n} - n | 0 \le n \le i-1 \})$, for $\ 0 \le i \le s$
and $j \ge \gr + i.$
\end{lemma}

\begin{proof} We use induction on $i$. For $i = 0$ the result is trivial.
Assume that $0\leq i\leq s-1$. We need to show that the inequality
holds for $i+1$. Let $j\geq \gr+i+1$. By (\ref{eqn2}), ${\mathfrak
{a}}_iI^j\cap a_{i+1}I^j=a_{i+1}[({\mathfrak
{a}}_iI^j:a_{i+1})\cap I^j]\subseteq a_{i+1}[({\mathfrak
{a}}_i:a_{i+1})\cap I^j]=a_{i+1}{\mathfrak {a}}_iI^{j-1}\subseteq
{\mathfrak {a}}_iI^j\cap a_{i+1}I^j$. Hence we obtain an exact
sequence
\begin{eqnarray}0\rightarrow a_{i+1}{\mathfrak {a}}_iI^{j-1}\rightarrow {\mathfrak
 {a}}_iI^j
\oplus a_{i+1}I^j\rightarrow {\mathfrak {a}}_{i+1}I^j\rightarrow
0. \label{eqn0}\end{eqnarray}
 On the other hand, by (\ref{eqn2}) for
$i=0$, $[0:a_{i+1}]\cap {\mathfrak {a}}_iI^{j-1}\subseteq
[0:a_{i+1}]\cap I^j=0$, and therefore $a_{i+1}{\mathfrak
{a}}_iI^{j-1}\cong {\mathfrak {a}}_iI^{j-1}$, $a_{i+1}I^j\cong
I^j$. The conclusion follows from (\ref{eqn0}) and the induction
hypothesis.
\end{proof}

\begin{lemma} \label{l2}{\rm \cite[2.5]{laura}}
Let $R$ be a local Cohen-Macaulay ring with infinite residue
field, and let $I$ be an $R$-ideal. Let $J$ be a reduction of $I$
with basis $a_1,\dots,a_s$ satisfying {\rm (\ref{eqn21})}. Then,
$\depth [\G/(a_1^*, \ldots, a_i^*)]_j \ge \min (\{ \depth R/I^n +
n - j-1 | j-i+1 \le n \le j+1 \}\cup \{\depth R/I^{j-i}-i+1\})$,
whenever $0 \le i \le s$ and $j \ge \gr+i+1$.
\end{lemma}

The goal of this section is to prove the following theorem.

\begin{thm} \label{reg}
Let $R$ be a local Cohen-Macaulay ring of dimension $d$ with
infinite residue field, let $I$ be an $R$-ideal, and let $J$ be a
reduction of $I$ generated by $s$ elements. Let $\G $ be the
associated graded ring of $I$, and $\gr = {\rm reg}\  \G$. Then,
\[ \depth \G \ge \min (\{ \depth R/I^j | 1 \le j \le \gr+1 \}\cup \{ \depth
R/I^j+j-\gr | 2+\gr \le j \le s+\gr \} ). \]
\end{thm}

To prove Theorem \ref{reg}, we will apply the methods of
\cite{laura}. We first need some preliminary notation and lemmas.

Let $J$ be a reduction of $I$ with basis $a_1, \ldots, a_s$
satisfying the conclusions of Lemma \ref{filter-regular}. If $s >
0$, for $0\le i\le s$ consider the graded $\G$-modules:
\[ M_{(i)} = [\G/(a_1^*, \ldots, a_i^*)]_{\ge \gr + i + 1} =
\sfrac{\G_+^{\gr+i+1}}{(a_1^*, \ldots, a_i^*)\G_+^{\gr+i}}
\]
\[ N_{(i)} = \sfrac{\G_+^{\gr+i}}{a_i^*\G_+^{\gr+i} + (a_1^*,
\ldots, a_{i-1}^*)\G_+^{\gr+i-1}}. \] Then, $[N_{(i)}]_{\ge
\gr+i+1} = M_{(i)}$
 and
$[N_{(i)}]_{\gr+i}=[\G/(a_1^*,\dots,a_{i-1}^*)]_{\gr+i}$.
 Hence, for $0\le i\le s$ we have the exact sequences
\begin{eqnarray}
 0 \rightarrow M_{(i)} \rightarrow N_{(i)} \rightarrow [\G/(a_1^*, \ldots,
a_{i-1}^*)]_{\gr+i} \rightarrow 0. \label{eqn3}\end{eqnarray}
Furthermore, if $0\leq i\leq s-1$, then
$N_{(i+1)}=M_{(i)}/a_{i+1}^*M_{(i)}$ and by (\ref{eqn21}) we have
that $0:_{M_{(i)}}(a_{i+1}^*)=0$. Thus, in the range $0\leq i\leq
s-1$ we have exact sequences \begin{eqnarray} 0\rightarrow
M_{(i)}(-1) \rightarrow M_{(i)}\rightarrow N_{(i+1)}\rightarrow 0.
\label{eqn33} \end{eqnarray}
 Notice that $M_{(s)}=0$, since
$I^{\gr +s+1}=JI^{\gr+s}$.

Let $\lambda=  \min (\{ \depth R/I^j | 1 \le j \le \gr+1 \}\cup \{
\depth R/I^j+j-\gr | 2+\gr \le j \le s+\gr \} ). $\newline Recall
that our goal is to show that $\depth \G\ge \lambda.$ The next
lemma gives an estimate of $\depth M_{(i)}$. In particular we show
that $\depth M_{(i)}\ge \lambda-i-1$.

\begin{lemma} \label{mi}
In addition to the assumptions of Theorem \ref{reg}, assume that
$s>0$. Let $M_{(i)}$ be defined as above. Then,
\[ \depth M_{(i)} \ge \min(\{ d-i, \depth R/I^{\gr + 1}-i + 1\}\cup \{ \depth R/I^{j} + j
-\gr-i-1 | 2+\gr \le j \le s+\gr \}). \]
\end{lemma}

\begin{proof} We use decreasing induction on $i$. For $i =s$, the assertion is true
since $M_{(s)} = 0$. Suppose that $0 \le i \le s-1$. Consider the
 exact sequence (\ref{eqn3})
\[ 0 \rightarrow M_{(i+1)} \rightarrow N_{(i+1)} \rightarrow [\G/(a_1^*, \ldots,
a_i^*)]_{\gr+i+1} \rightarrow 0. \] It follows from Lemma \ref{l2}
and the induction hypothesis that $\depth N_{(i+1)}\ge
\min(\{d-i-1, \depth R/I^{\gr + 1} -i\}\cup\{ \depth
R/I^j+j-\gr-i-2 | 2+\gr \le j \le s+\gr \}\cup\{\depth
R/I^j+j-\gr-i-2 | 2+\gr \le j \le \gr + i+2 \}).$ If $i\le s-2$,
then $\depth N_{(i+1)}\ge \min(\{d-i-1, \depth R/I^{\gr + 1}
-i\}\cup\{ \depth R/I^j+j-\gr-i-2 | 2+\gr \le j \le s+\gr \}).$ If
$i= s-1$, then by Lemma \ref{l1} we have that $\depth
R/I^{\gr+s+1} = \depth R/JI^{\gr+s} \ge \min (\{ d-s\}\cup\{
\depth R/I^{\gr+s-n}-n | 0 \le n \le s-1 \}).$ It follows that
$\depth N_{(s)}\ge \min(\{d-s, \depth R/I^{\gr + 1} -s+1\}\cup\{
\depth R/I^j+j-\gr-s-1 | 2+\gr \le j \le s+\gr \}).$ In any case
$\depth N_{(i+1)}\ge \min(\{d-i-1, \depth R/I^{\gr + 1} -i\}\cup\{
\depth R/I^j+j-\gr-i-2 | 2+\gr \le j \le s+\gr \})$, and since
$\depth M_{(i)} = \depth N_{(i+1)} + 1$, the conclusion follows.
\end{proof}

Let $S$ be a homogeneous Noetherian ring with $S_0$ local and
homogeneous maximal ideal $\mathfrak {M}$, let $H^{\bullet}(-)$
denote local cohomology with support in $\mathfrak {M}$.
\newline For a graded $S$-module $N$ and an integer $j$ we put
$a_j(N)=\max\{n \mid [H^j(N)]_n\neq 0\}$.\newline The following
lemma is well known, but we recall it for convenience.

\begin{lemma}\label{lemma6}{\rm \cite[2.6]{laura}}
Let $0\rightarrow A \rightarrow B \rightarrow
C\rightarrow 0$ be an exact sequence of graded $S$-modules, let
$n$ and $j$ be integers.
\begin{enumerate}
\item[(a)]
If $a_j(A)\leq n$ and $a_j(C)\leq n$, then $a_j(B)\leq n$.
\item[(b)]
\begin{enumerate}
\item[(i)] If $H^j(A)=0$, then $a_j(C)\geq a_j(B)$.

\item[(ii)]If $H^j(B)=0$, then $a_{j+1}(A)\geq a_j(C)$.
\item[(iii)] If $H^j(C)=0$, then $a_{j+1}(B)\geq a_{j+1}(A)$.
\end{enumerate}
\end{enumerate}
\end{lemma}

\begin{lemma} \label{laura27}
In addition to the assumptions of Theorem \ref{reg}, assume that
$s>0$. Let $M_{(i)}$ and $\lambda$ be defined as above. Then,
\begin{enumerate}
\item $a_j( M_{(i)}) \le \gr+i$ for any integer $j$ and $0\le i\le
s$.
\item $\depth  M_{(i)} \ge \lambda-i-1$ and if $\depth M_{(i)} = \lambda - i-1$ then
$a_{\lambda-i-1}(M_{(i)}) = \gr+i$.
\end{enumerate}
\end{lemma}

\begin{proof}
{\bf (1)} We prove the claim by decreasing induction on $i$. For
$i =s$ the assertion is trivial, since $M_{(s)}=0$. Suppose that
$0 \le i \le s-1$ and that $a_j( M_{(i+1)}) \le \gr+i+1$ for any
integer $j$. Consider the exact sequence (\ref{eqn3})
\[ 0 \rightarrow M_{(i+1)} \rightarrow N_{(i+1)} \rightarrow [\G/(a_1^*, \ldots,
a_i^*)]_{\gr+i+1} \rightarrow 0. \] By \cite[2.2]{GH}, for any
integer $j$, $H^j([\G/(a_1^*, \ldots, a_i^*)]_{\gr+i+1})$ is
concentrated in degree $\gr+i+1$. Hence, Lemma~\ref{lemma6}
$(${\rm a}$)$ implies that $ a_j(N_{(i+1)}) \le \gr+i+1$ for any
$j$. Applying the local cohomology functor to the sequence
(\ref{eqn33}) $$ 0\rightarrow M_{(i)}(-1) \rightarrow
M_{(i)}\rightarrow N_{(i+1)}\rightarrow 0
$$
 it follows that for any $j$, $[H^j(M_{(i)})]_n=0$ whenever $n>\gr+i$. Hence
$a_j(M_{(i)})\leq \gr+i$ and the proof of (a) is completed.

{\bf (2)} It follows from Lemma \ref{mi} that $\depth M_{(i)} \ge
\lambda-i-1$. To prove the last assertion, we again use decreasing
induction on $i$. For $i = s$, the assertion is vacuous. Suppose
that $0 \le i \le s-1$, and that $\depth M_{(i)}=\lambda-i-1$. It
follows from (\ref{eqn33}) that $\depth N_{(i+1)} = \lambda-i-2$,
and so $H^{\lambda-i-2}(N_{(i+1)}) \not= 0$. Applying the local
cohomology functor to (\ref{eqn3}) we obtain the exact sequence
\[ \dots \rightarrow  H^{\lambda-i-2}(M_{(i+1)}) \rightarrow H^{\lambda-i-2}(N_{(i+1)}) \rightarrow
H^{\lambda-i-2}\big([\G/(a_1^*, \ldots, a_i^*)]_{\gr+i+1}\big)
\rightarrow \dots. \] If $\depth M_{(i+1)} > \lambda-i-2$, then
$H^{\lambda-i-2}(N_{(i+1)}) \hookrightarrow
H^{\lambda-i-2}([\G/(a_1^*, \ldots, a_i^*)]_{\gr+i+1})$, and so
$a_{\lambda-i-2}(N_{(i+1)}) = \gr+i+1$. If $\depth M_{(i+1)} =
\lambda-i-2$, then by induction hypothesis we have that
$a_{\lambda-i-2}(M_{(i+1)}) = \gr+i+1$. Again, we consider the
exact sequence (\ref{eqn3})
\[ 0 \rightarrow M_{(i+1)} \rightarrow N_{(i+1)} \rightarrow [\G/(a_1^*, \ldots,
a_i^*)]_{\gr+i+1} \rightarrow 0. \] It follows from Lemma \ref{l2}
(and Lemma \ref{l1} when $i = s-1$) that $ \depth [\G/(a_1^*,
\ldots, a_i^*)]_{\gr+i+1} \ge \lambda-i-2$. Thus, applying
Lemma~\ref{lemma6} $(${\rm b}$)$  $(${\rm iii}$)$ with
$j=\lambda-i-3$ to the exact sequence above
we get that
$a_{\lambda-i-2}(N_{(i+1)})\geq \gr+i+1$. In any case, we have
that $a_{\lambda-i-2}(N_{(i+1)}) \ge \gr+i+1$. Since depth
$M_{(i)}=\lambda-i-1$, applying Lemma~\ref{lemma6} $(${\rm b}$)$
$(${\rm ii}$)$ to the sequence (\ref{eqn33}) we conclude that
$a_{\lambda-i-1}(M_{(i)})\geq \gr+i$.
\end{proof}

We are now ready to prove Theorem \ref{reg}.

\begin{proofc} We need to show that $\depth \G\ge \lambda$. Let $J$ be a reduction of $I$ with basis
$ a_1, \ldots, a_s $
satisfying the conclusions of Lemma \ref{filter-regular}. If $s=
0$, then $\G = R/I \oplus I/I^2 \oplus \ldots \oplus
I^{r-1}/I^r\oplus I^r, $ where $r = r_J(I)$ is the reduction
number of $I$ with respect to $J$. Hence $\depth \G=\min\{\depth
R/I^j\mid 1 \le j\le r\}$. The result follows from the fact that
$\gr \ge r$.\newline Suppose now that $s > 0$. From the definition
of $M_{(0)}$, we have the exact sequence
\begin{eqnarray}
 0  \rightarrow  M_{(0)} \rightarrow  \G  \rightarrow
\oplus_{n=0}^\gr I^n/I^{n+1}  \rightarrow  0. \label{eq5}
\end{eqnarray}
Let $C = \oplus_{n=0}^\gr I^n/I^{n+1}$. Since $\depth C \ge
\lambda$, and $\depth M_{(0)} \ge \lambda - 1$ by  Lemma
\ref{laura27}, it follows that $\depth \G \ge \lambda-1$. Applying
local cohomology to (\ref{eq5}) we see that $
 H^{\lambda-1}(M_{(0)})\cong
H^{\lambda-1}(\G)$. Furthermore, by Lemma \ref{laura27} and by
Lemma~\ref{lemma6} $(${\rm a}$)$, we have that $a_j(\G)\le \gr$
for any integer $j$. If $\depth \G=\lambda -1 $, then $\depth
M_{(0)} = \lambda-1$, and so $a_{\lambda-1}(M_{(0)}) = \gr$ by
Lemma \ref{laura27}. On the other hand, since $\lambda -1 < d$, by
\cite[2.9]{laura} we have that $a_{\lambda-1}(\G) < \gr$, a
contradiction. Hence, $\depth \G\ge\lambda$.
\end{proofc}

\end{document}